\documentclass[11pt]{article}
\usepackage{times,mathptm}
\usepackage{amsfonts,amscd,amssymb}
\input diagrams
\newtheorem{definition}{Definition}[section]
\newtheorem{lemma}[definition]{Lemma}
\newtheorem{proposition}[definition]{Proposition}
\newtheorem{corollary}[definition]{Corollary}
\newtheorem{remark}[definition]{Remark}
\newtheorem{remarks}[definition]{Remarks}
\newtheorem{theorem}[definition]{Theorem}
\newtheorem{example}[definition]{Example}
\newtheorem{examples}[definition]{Examples}
\newcommand{\thlabel}[1]{\label{th:#1}}
\newcommand{\thref}[1]{Theorem~\ref{th:#1}}
\newcommand{\selabel}[1]{\label{se:#1}}
\newcommand{\seref}[1]{Section~\ref{se:#1}}
\newcommand{\lelabel}[1]{\label{le:#1}}

\newcommand{\prlabel}[1]{\label{pr:#1}}
\newcommand{\prref}[1]{Proposition~\ref{pr:#1}}
\newcommand{\colabel}[1]{\label{co:#1}}
\newcommand{\coref}[1]{Corollary~\ref{co:#1}}
\newcommand{\relabel}[1]{\label{re:#1}}

\newcommand{\exlabel}[1]{\label{ex:#1}}
\newcommand{\exref}[1]{Example~\ref{ex:#1}}
\newcommand{\delabel}[1]{\label{de:#1}}
\newcommand{\deref}[1]{Definition~\ref{de:#1}}
\newcommand{\eqlabel}[1]{\label{eq:#1}}
\newcommand{\eqref}[1]{(\ref{eq:#1})}

\newenvironment{proof}{{\it Proof.}}{\hfill $ \square $ \vskip 4mm}
\newcommand{\Hom}{\rm{Hom}\,}

\newcommand{\End}{\rm{End}\,}

\newcommand{\Cc}{\mathcal{C}}
\newcommand{\Dd}{\mathcal{D}}
\newcommand{\Ee}{\mathcal{E}}
\newcommand{\Ff}{\mathcal{F}}
\newcommand{\Hh}{\mathcal{H}}
\newcommand{\Mm}{\mathcal{M}}
\newcommand{\Pp}{\mathcal{P}}

\def\text#1{\mbox{{\rm #1}}}

\def\ol{\overline}

\def\dul#1{\underline{\underline{#1}}}
\def\Nat{\dul{\rm Nat}}

\def\ot{\otimes}

\def\doublerightleft#1#2{{\lower.2ex\vbox{
\hbox{${\smash{\mathop{\longrightarrow}\limits^{#1}}}$}\vspace*{-4mm}
\hbox{${\smash{\mathop{\longleftarrow}\limits_{#2}}}$}}}}

\begin{document}
\title{Maschke functors, semisimple functors
and separable functors of the second kind. Applications
\thanks{Research supported by the bilateral project
``Hopf Algebras in Algebra, Topology, Geometry and Physics" of the Flemish
and
Romanian governments.}}
\author{
S. Caenepeel
\\ Faculty of Applied Sciences\\
Free University of Brussels, VUB\\ B-1050 Brussels, Belgium\and
G. Militaru\\
Faculty of Mathematics\\ University of Bucharest\\
RO-70109 Bucharest 1, Romania}
\date{}
\maketitle

\begin{abstract}
We introduce separable functors of the second kind
(or $H$-separable functors) and $H$-Maschke functors.
$H$-separable functors are generalizations of separable functors.
Various necessary and sufficient conditions for a functor to be
$H$-separable or $H$-Maschke, in terms of generalized
(co)Casimir elements (integrals, in the case of Hopf algebras),
are given. An $H$-separable functor is always $H$-Maschke,
but the converse holds in particular situations.
A special role will be played by Frobenius functors and
their relations to $H$-separability.
Our concepts are applied to modules,
comodules, entwined modules, quantum Yetter-Drinfeld modules,
relative Hopf modules.
\end{abstract}

\section{Introduction}\selabel{0}
One of the fundamental results in classical representation theory
is {\sl Maschke's Theorem}, stating that a finite group algebra
$kG$ over a field $k$ is {\sl semisimple} if and only if the characteristic
of $k$ does not divide the order of $G$. Several generalizations of
this result have appeared in the literature. To illustrate that there
is a subtle difference, let us look more carefully at one of the
earliest generalization, where the ground field $k$ is replaced by
a commutative ring $k$. Algebras over a commutative ring are rarely
semisimple, and one arrives at the following result:
a finite group algebra
$kG$ over a field $k$ is {\sl separable} if and only if the characteristic
of $k$ does not divide the order of $G$. The interesting thing is
that, over a field $k$, a separable finite dimensional algebra is
semisimple, but not conversely: it suffices to look at a purely
inseparable field extension. A consequence of the two versions of
Maschke's Theorem is then that, for a finite group algebra (and,
more generally, for a finite dimensional Hopf algebra)
over a field, separability and
semisimplicity are equivalent.\\
An elegant categorical definition of separability has been proposed
by N\v ast\v asescu et al. in \cite{NastasescuVV89}. A functor
$F$ is called separable if and only if the natural transformation
$\Ff$ induced by $F$ is split by a natural transformation $\Pp$.
It is a proper
generalization of the notion of separable algebra, in the sense that
a
$k$-algebra $A$ is separable if and only if the restriction of
scalars functor $F: {\cal M}_A \to {\cal M}_k$ is
separable \cite[Prop. 1.3]{NastasescuVV89}.
Moreover, a separable functor $F$ between two abelian
categories satisfies the following version of Maschke's Theorem:
an exact sequence, that becomes split after applying $F$ is itself
split. If we apply this property to the restriction of scalars
functor in the case of an algebra $A$ over a field $k$, then we
easily deduce that this algebra is semisimple. We also point out that
many of the recent Maschke-type Theorems (see e.g.
\cite{CaenepeelMZ97a}, \cite{CaenepeelIMZ98}, \cite{book},...)
come down to proving that a
certain functor is separable.\\
Now consider a separable algebra $A$ over a field. What are
its properties that distinguish it from a semisimple algebra?
The answer is the following: $\Pp$ allows
to deform a $k$-linear splitting map $f$ between two $A$-modules in
such a way that it becomes $A$-linear; this can also be done in the
semisimple case, but in the separable case the deformation is natural
in $f$!\\
In \seref{3}, we will propose categorical properties of functors that,
when applied to the restriction of scalars functor in the case of an algebra
over a field $k$, are equivalent to semisimplicity of the algebra.
The starting point is the following: an algebra $A$ over a field is
semisimple
if and only if every $A$-module is projective, if and only if every
$A$-module is injective. We will say that a functor $F:\Cc \to \Dd$ is a
Maschke (resp. dual Maschke) functor if every object in $\Cc$ is relative
injective (resp. projective). A functor $F$ between abelian categories
is called semisimple if and only if an exact sequence that becomes split
after applying $F$ is itself split. The three notions are equivalent for
a functor reflecting monics and epics (see \prref{3.7}).\\
In \seref{2}, we introduce another generalization of separable functors;
consider an exact sequence of graded modules over a $G$-graded
$k$-algebra $A$. Suppose that the sequence is split after we forget
the $A$-action and the $G$-grading; separability of the functor forgetting
action and grading would imply that the sequence is split as a sequence
of graded $A$-modules. When can we conclude that the sequence is at
least split as a sequence of $A$-modules? Or consider the following:
an exact sequence of $A$-modules, with $A$ a $k$-algebra, which is
split as a sequence of $k$-modules. Is it split as a sequence of
$B$-modules, where $B$ is a given subalgebra of $A$.\\
This leads us to the following: let $F:\ \Cc\to \Dd$ and
$H:\ \Cc\to \Ee$ be functors. We $F$ an $H$-separable functor,
or separable functor of the second kind, if the
natural transformation $\Hh$ induced by $H$ factorizes as a
natural transformation trough $\Ff$ induced by $F$:
$$\Hh=\Pp\circ \Ff, $$
for a natural transformation $\Pp$.
If $H$ is the identity functor, then we recover the separable functors
of \cite{NastasescuVV89}. Most properties of separable functors
(Maschke's Theorem, Rafael's Theorem, the Frobenius-Rafael Theorem) can be
generalized to $H$-separable functors, we discuss this in \seref{2}.
Also the notion of Maschke functor, dual Maschke functor, and
semisimple functor can be generalized in the same spirit; in fact,
we decided to present at once the general theory of
(dual) $H$-Maschke functors in \seref{3}.\\
In \seref{4}, we present some examples and applications, we look at
the categories of modules, comodules, entwined modules, Hopf modules
and relative Hopf modules. We present a structure Theorem for
injective objects in the category of entwined modules, that
arose from noncommutative geometry \cite{BrzezinskiM}.
A separable
functor is always Maschke (and dual Maschke), and in some particular
cases we have the converse property. As a first example, we have
modules over a group algebra or a Hopf algebra. The fact that Maschke
implies separability comes from the fact that the separability of
a Hopf algebra can be described in terms of integrals in the
Hopf algebra. A similar phenomenon appears if we look at relative Hopf
modules: if the functor forgetting action and coaction is $H$-Maschke
($H$ is the functor forgetting $A$-action), then it is also
$H$-separable. Both conditions (\thref{4.10})
are equivalent to the fact that
there exists a total integral in the sense of Doi \cite{Doi85}.
This gives another motivation for introducing the $H$-separability
concept.

\section{Preliminaries}\selabel{1}
Let $k$ be a commutative ring. For a $k$-coalgebra $C$, we use
the Sweedler-Heyneman notation for the comultiplication $\Delta_C$:
$$\Delta_C(c)=c_{(1)}\ot c_{(2)}$$
(summation is implicitely understood). $\varepsilon_C$ will denote
the counit of $C$. If $C$ coacts from the right on a $k$-module $M$,
that is, $M$ is a right $C$-comodule, then we write $\rho_M$ for the
structure map, and use the following notation:
$$\rho_M(m)=m_{[0]}\ot m_{[1]} \in M\ot C$$
${\cal M}^C$
will be the category of right $C$-comodules and
$C$-colinear maps. A right $C$-comodule
$M$ is called {\it relative injective}  if for any
$k$-split monomorphism $i:U\to V$ in ${\cal M}^C$ and for any
$C$-colinear map $f:U\to M$, there exists a $C$-colinear map
$g:V\to M$ such that $g\circ i=f$. This is equivalent
to the fact that $\rho_M: M\to M\ot C$ splits in ${\cal M}^C$,
i.e. there exists
a $C$-colinear map $\lambda_M :M\ot C\to M$ such that
$\lambda_M \circ \rho_M= {\rm I}_M$. Of course, if $k$ is a
field, $M$ is relative injective if and only if it is an injective object
in ${\cal M}^C$. Relative projective modules over a $k$-algebra
$A$ are defined dually. \\
A (right-right) entwining structure (\cite{BrzezinskiM})
is a triple $(A,C,\psi)$, where
$A$ is a $k$-algebra, $C$ is a $k$-coalgebra, and $\psi:\
C\ot A\to A\ot C$ is a $k$-linear map satisfying the conditions
\begin{eqnarray}
&&  (ab)_{\psi}\ot c^{\psi}=  a_{\psi}b_{\Psi}\ot
c^{\psi\Psi}\eqlabel{1.1}\\
&&(1_A)_{\psi}\ot c^{\psi}=1_A\ot c\eqlabel{1.2}\\
&&  a_{\psi}\ot \Delta_C(c^{\psi})=
  a_{\psi\Psi}\ot c_{(1)}^{\Psi}\ot c_{(2)}^{\psi}
\eqlabel{1.3}\\
&&\varepsilon_C(c^{\psi})  a_{\psi} =
\varepsilon_C(c)a\eqlabel{1.4}
\end{eqnarray}
Here we used the sigma notation
$$\psi(c\ot a)=  a_{\psi}\ot c^{\psi}
=  a_{\Psi}\ot c^{\Psi}$$
Entwining structures where introduced with a motivation
coming from noncommutative geometry: one can generalize
the notion of principal bundles to a very general setting
in which the role
of coordinate functions on the base is played by a general
noncommutative algebra $A$, and the fibre of the principal bundle
by a coalgebra $C$, where $A$ and $C$ are related by a
map $\psi:\ A\ot C\to C\ot A$, called the entwining map.

An entwining module $M$ is at the same time a right $A$-module
and a right $C$-comodule such that the following compatibility
relation holds between the action and coaction:
$$\rho(ma)=m_{[0]}a_{\psi}\ot m_{[0]}^{\psi}$$
$\Mm(\psi)_A^C$ is the category of entwined modules and $A$-linear
$C$-colinear maps. The forgetful functors
$$F:\ \Mm(\psi)_A^C\to \Mm_A~~{\rm and}~~H:\ \Mm(\psi)_A^C\to \Mm^C$$
have respectively a right and a left adjoint
(\cite{Brzezinski1})
$$G=\bullet\ot C~~{\rm and}~~K=\bullet \ot A.$$
A Doi-Koppinen datum $(H,A,C)$ consists of a bialgebra $H$, a right
$H$-comodule algebra $A$ and a right $H$-module coalgebra $C$. Associated to
it is an entwined structure $(A,C,\psi)$, with
$$\psi(c\ot a)=a_{[0]}\ot ca_{[1]}$$
The following special cases will be of interest to us:\\
1) $C=H$, where $H$ is a Hopf algebra. In this case
$\Mm(\psi)_A^C=\Mm_A^H$, the category of relative Hopf modules.\\
2) $A=H$, where again $H$ is a Hopf algebra. Now
$\Mm(\psi)_A^C=\Mm_H^C$, the category of Doi's $[H,C]$-modules.

\section{$H$-separable functors}\selabel{2}
Let $F:\ \Cc\to\Dd$ and $H:\ \Cc\to\Ee$ be covariant functors.
We then have functors
$$\Hom_{\Cc}(\bullet,\bullet),~\Hom_{\Dd}(F,F),~~\Hom_{\Ee}(H,H):\
\Cc^{\rm op}\times \Cc\to\dul{\rm Sets}$$
and natural transformations
$$\Ff:\ \Hom_{\Cc}(\bullet,\bullet)\to \Hom_{\Dd}(F,F)~~;~~
\Hh:\ \Hom_{\Cc}(\bullet,\bullet)\to \Hom_{\Ee}(H,H)$$
given by
$$\Ff_{C,C'}(f)=F(f)~~;~~\Hh_{C,C'}(f)=H(f)$$
for $f:\ C\to C'$ in $\Cc$.

\begin{definition}\delabel{2.1}\rm
The functor $F$ is called $H$-separable if there exists a natural
transformation
$$\Pp:\ \Hom_{\Dd}(F,F)\to \Hom_{\Ee}(H,H)$$
such that
\begin{equation}\eqlabel{2.1.1}
\Pp\circ \Ff=\Hh
\end{equation}
that is, $\Hh$ factors through $\Ff$ as a natural transformation,
and we have a commutative diagram
$$\begin{diagram}
\Hom_{\Cc}(\bullet,\bullet)&\rTo^{\Ff}&\Hom_{\Dd}(F,F)\\
\dTo_{\Hh}&\SW_{\Pp}\\
\Hom_{\Ee}(H,H)
\end{diagram}$$
\end{definition}

\begin{remarks}\relabel{2.2}\rm
1) $F$ is $1_{\Cc}$-separable if and only if $F$ is separable in
the sense of \cite{NastasescuVV89}. Indeed, the functor $F$ is
separable if and only if there exists a natural transformation $\Pp$
such that $\Pp\circ \Ff=1_{\Cc}$ (see \cite{book}). We refer to
\cite{book} for a detailed study of separable functors.
A  finite extension of commutative fields
$k\subset K$ is separable in the classical sense if
and only if the forgetful
functor $F :{\cal M}_K \to {\cal M}_k$ is separable.
For the reader convenience we show how the above
natural transformation ${\cal P}$ is constructed:
let $K/k$ be a finite
separable extension, $\alpha \in K$ be a primitive
element (i.e. $K=k(\alpha)$) and $p\in K[X]$,
$p(X) = X^n - \sum_{i=0}^{n-1} c_i X^i$ be the
minimal polynomial of $\alpha$. Then the
natural transformation ${\cal P}$ is constructed
as follows: for $M$, $N$ two $K$-vector space
we define
$$
{\cal P}_{M, N} :\Hom_k (M, N) \to \Hom_K (M, N), \quad
{\cal P}_{M, N} (f) (m): =
p'(\alpha)^{-1} \sum_{i=0}^{n-1}
\alpha^{-i-1}(\sum_{j=0}^i c_j \alpha^j) f(\alpha^{i} m)
$$
for any $f\in \Hom_k (M, N)$ and $m\in M$. Then
${\cal P}$ is a natural transformation that splits
${\cal F}$. This is a one of remarkable property of
classical separable fields extension $K/k$: any
$k$-linear map $f$ between two $K$-vector spaces
can be deformated, using the above
formula, until it becomes a $K$-linear map.\\
As we will see below, most properties of separable functors
can be generalized to $H$-separable functors.\\
2) The fact that $\Pp$ is natural means the following condition:
for
$$u:\ X\to Y,~v:\ Z\to T~{\rm in}~\Cc~~{\rm and}~~
h:\ F(Y)\to F(Z)~{\rm in}~\Dd,$$
we have
\begin{equation}\eqlabel{2.2.1}
\Pp_{X,T}(F(v)\circ h\circ F(u))=
H(v)\circ \Pp_{Y,Z}(h)\circ H(u)
\end{equation}
\eqref{2.1.1} can be rewritten as
\begin{equation}\eqlabel{2.2.2}
\Pp_{C,C'}(F(f))=H(f)
\end{equation}
for any $f:\ C\to C'$ in $\Cc$.
\end{remarks}

\begin{proposition}\prlabel{2.3}
Consider functors
$$\Cc\rTo^{F}\Dd\rTo^{F_1}\Dd_1~~{\rm and}~~\Cc\rTo^{H}\Ee$$
1) If $F_1\circ F$ is $H$-separable, then $F$ is $H$-separable.\\
2) If $F$ is $H$-separable, and $F_1$ is separable, then
$F_1\circ F$ is $H$-separable.
\end{proposition}

\begin{proof}
Obvious.
\end{proof}

\begin{proposition}\prlabel{2.4}
Let $F$ be an $H$-separable functor. If $f:\ C\to C'$ in $\Cc$ is
such that $F(f)$ has a left, right, or two-sided inverse in $\Dd$,
then $H(f)$ has a left, right, or two-sided inverse in $\Ee$.
\end{proposition}

\begin{proof}
Let $g$ be a left inverse of $F(f)$. Using \eqref{2.2.1} and
\eqref{2.2.2}, we find
\begin{eqnarray*}
\Pp_{C,C'}(g)\circ H(f)&=&
\Pp_{C,C'}(g\circ F(f))=\Pp_{C,C'}(I_{F(C)})\\
&=& \Pp_{C,C'}(F(I_C))=H(I_C)=I_{H(C)}
\end{eqnarray*}
The proof for right and two-sided inverses is similar.
\end{proof}

\begin{corollary}\colabel{2.5}
{\bf (Maschke's Theorem for $H$-separable functors)}
Let $\Cc$, $\Dd$ and $\Ee$ be abelian categories, and assume that
$F:\ \Cc\to\Dd$ is $H$-separable. An exact sequence in $\Cc$ that
becomes split after we apply the functor $F$, also becomes split
after we apply the functor $H$.
\end{corollary}

Recall that Rafael's Theorem (see \cite{Rafael90}) gives an easy
criterion for the separability a functor that has a left or right adjoint.
We will now generalize Rafael's Theorem to $H$-separable functors.
First, we recall the following well-known result from category theory.
For completeness sake, we include a brief sketch of proof, based on
the well-konwn property that $(F,G)$ a pair of adjoint functors
between the categories $\Cc$ and $\Dd$ if and only if
there exist two
natural transformations $\eta:\ 1_\Cc\to GF$
and $\varepsilon:\ FG\to 1_\Dd$, called the unit and
counit of the adjunction, such that
\begin{equation}\eqlabel{2.5.1}
G(\varepsilon_D)\circ\eta_{G(D)}= I_{G(D)}~~{\rm and}~~
\varepsilon_{F(C)}\circ F(\eta_C)= I_{F(C)}
\end{equation}
for all $C\in \Cc$ and $D\in \Dd$.

\begin{proposition}\prlabel{2.6}
Let $G:\ \Dd\to \Cc$ be a right adjoint of $F:\ \Cc\to \Dd$,
and consider functors $H:\ \Cc\to \Ee$ and $K:\ \Dd\to \Ee$.
Then we have isomorphisms
\begin{eqnarray*}
\Nat(HGF,H)&\cong& \Nat\Bigl(\Hom_{\Dd}(F,F),\Hom_{\Ee}(H,H)\Bigr)\\
\Nat(K,KFG)&\cong&\Nat\Bigl(\Hom_{\Cc}(G,G),\Hom_{\Ee}(K,K)\Bigr)
\end{eqnarray*}
\end{proposition}

\begin{proof}
For $\nu:\ HGF\to H$, we define
$${\cal P}=\alpha(\nu):\ \Hom_{\Dd}(F,F)\to \Hom_{\Ee}(H,H)$$
as follows: for $g:\ F(C)\to F(C')$ in $\Dd$, we put
$$\Pp_{C,C'}(g)=\nu_{C'}\circ HG(g)\circ H(\eta_C)$$
Conversely, given $\Pp:\ \Hom_{\Dd}(F,F)\to \Hom_{\Ee}(H,H)$,
we define $\alpha^{-1}(\Pp):\ HGF\to H$ by
$$\nu_C=\Pp_{GF(C),C}(\varepsilon_{F(C)})$$
for any $C\in \Cc$.
\end{proof}

\begin{theorem} \thlabel{2.7}
{\bf (Rafael's Theorem for $H$-separability)}
Let $G:\ \Dd\to \Cc$ be a right adjoint of $F:\ \Cc\to \Dd$,
and consider functors $H:\ \Cc\to \Ee$ and $K:\ \Dd\to \Ee$.
Then:\\
1) $F$ is $H$-separable if and only if there exists a natural
transformation $\nu:\ HGF\to H$ such that
\begin{equation}\eqlabel{2.7.1}
\nu_C\circ H(\eta_C)=I_{H(C)}
\end{equation}
for any $C\in \Cc$.\\
2) $G$ is $K$-separable if and only if there exists a natural
transformation $\zeta:\ K\to KFG$ such that
\begin{equation}\eqlabel{2.7.2}
K(\varepsilon_D)\circ \zeta_D=I_{K(D)}
\end{equation}
for any $D\in \Dd$.
\end{theorem}

\begin{proof} We only prove the first statement; the proof of the
second one is similar. We use the notation introduced in the proof
of \prref{2.6}.
Assume that $F$ is $H$-separable, and put $\nu=\alpha^{-1}(\Pp)$.
Then we compute
\begin{eqnarray*}
\nu_C\circ H(\eta_C)&=& \Pp_{GF(C),C}(\varepsilon_{F(C)})\circ H(\eta_C)\\
{\rm \eqref{2.2.1}}~~~~&=&
\Pp_{C,C}(\varepsilon_{F(C)}\circ F(\eta_C))\\
{\rm \eqref{2.5.1}}~~~~&=&\Pp_{C,C}(I_{F(C)})=\Pp_{C,C}(F(I_C))=
H(I_C)=I_{H(C)}
\end{eqnarray*}
Conversely, assume that $\nu$ satisfies \eqref{2.7.1}, and take
$\Pp=\alpha({\nu})$. Using \eqref{2.2.1}, we find
$$\Pp(F(f))=\nu_{C'}\circ HGF(f)\circ H(\eta_C)=
H(f)\circ\nu_C\circ H(\eta_C)=H(f)$$
as needed.
\end{proof}

Recall \cite{CaenepeelMZ97b} that a functor $F$
is called Frobenius if
$F$ has a right adjoint $G$ that is also a left adjoint.
$(F,G)$ is then called a Frobenius pair. In
\cite{CaenepeelK01}, a Rafael-type criterion for the separability
of a Frobenius functor is given. We will now generalize this to
$H$-separability. First, we need the following standard fact
from category theory.

\begin{proposition}\prlabel{2.8}
Let $G$ be a left adjoint of the functor $F:\ \Cc\to \Dd$,
and $H:\ \Cc\to \Ee$ a functor. Then we have an isomorphism
$$\Nat(HGF,H)\cong \Nat(HG, HG)$$
\end{proposition}

\begin{proof} (sketch)
Let
$$\mu:\ GF\to 1_{\Cc}~~{\rm and}~~\chi:\ 1_{\Dd}\to FG$$
be the counit and unit of the adjunction $(G,F)$.
For $\nu:\ HGF\to H$,
we define $\beta=\beta_{\nu}: HG\to HG$ by
$$\beta_D=\nu_{G(D)}\circ HG(\chi_D):\ HG(D)\to HG(D)$$
for every $D\in \Dd$. Conversely, given
$\beta:\ HG\to HG$, we define
$\nu=\nu_{\beta}:\ HGF\to H$ as follows:
$$\nu_C = H(\mu_C)\circ\beta_{F(C)}:\
HFG(C)\to H(C)$$
for every $C\in \Cc$.
\end{proof}

For a Frobenius pair of functors $(F,G)$, we will write
$$\chi:\ 1_{\Dd}\to FG ~~{\rm and}~~ \mu:\ GF\to 1_{\Cc}$$
be the unit and counit of the adjunction $(G,F)$ and
$$\eta:\ 1_\Cc\to GF ~~{\rm and}~~ \varepsilon:\ FG\to 1_\Dd$$
the unit and counit of the adjunction $(F, G)$.

\begin{proposition}\prlabel{2.9}
Let $(F,G)$ be a Frobenius pair and $H:\ \Cc\to \Ee$
a functor. Then $F$ is $H$-separable if
and only if there exists a natural transformation
$\beta:\ HG\to HG$ such that
\begin{equation}\eqlabel{2.9.1}
H(\mu_C)\circ\beta_{F(C)}\circ H(\eta_C)=I_{H(C)}
\end{equation}
for all $C\in \Cc$.
\end{proposition}

\begin{proof}
First we will apply \thref{2.7} to the adjunction $(F, G)$:
we obtain that $F$ is $H$-separable if and only if there exists
a natural transformation $\nu:\ HGF\to H$ such that
\eqref{2.7.1} holds.

Now, we apply \prref{2.8} to the adjunction $(G,F)$ to obtain
the corresponding natural transformation $\beta =\beta_{\nu}$.
Furthermore, \eqref{2.7.1} holds for $\nu$ if and only if
\eqref{2.9.1} holds for $\beta =\beta_{\nu}$.
\end{proof}

\section{Relative injectivity and Maschke functors}\selabel{3}
\begin{definition}\delabel{3.1}\rm
Let $F :\ \Cc \to \Dd$ and $H:\ \Cc\to \Ee$ be
covariant functors. An object
$M \in \Cc$ is called {\em $F$-relative $H$-injective} if the
following condition is satisfied: for any $i:\ C\to C'$
in $\Cc$ with $F(i) :\ F(C)\to F(C')$ a split monic in $\Dd$,
and for every $f:\ C\to M$ in $\Cc$, there exists
$g :\ H(C')\to H(M)$ in $\Ee$ such that $H(f)=g\circ H(i)$, that is,
the following diagram commutes in $\Ee$:
\begin{equation}\eqlabel{diag1}
\begin{diagram}
 H(C)       &  \rTo^{H(i)}  &   H(C') \\
\dTo^{H(f)} &  \SW_g \\
H(M)        &          &
\end{diagram}
\end{equation}
$F$ is called an {\em $H$-Maschke functor} if any object of $\Cc$ is
$F$-relative $H$-injective.\\
An $F$-relative $1_{\Cc}$-injective object is also called an
$F$-relative injective object. A $1_{\Cc}$-Maschke functor is also
called a {\em Maschke functor}.\\
$P \in \Cc$ is called {\em $F$-relative $H$-projective} if
$P$ is $F^{\rm op}$-relative $H^{\rm op}$-injective, where
$F^{\rm op}:\ \Cc^{\rm op}\to \Dd^{\rm op}$ is the functor
opposite to $F$.\\
$F$ is called a {\em dual $H$-Maschke functor} if any object of $\Cc$ is
$F$-relative $H$-projective.
\end{definition}

\begin{examples}\exlabel{3.2}\rm
1) Every
object of $\Cc$ is $1_{\Cc}$-relative injective.\\
2) Let $A$ be an algebra over a field $k$, $\Dd$ the category
of $k$-vector spaces, and $\Cc={\cal M}_A$ the category of right
$A$-modules (or representations of $A$). The restrictions of scalars
functor $F: {\cal M}_A\to {\cal M}_k$ is exact,
and every monic (resp. epic) in ${\cal M}_k$
splits (resp. cosplits), and therefore an $A$-module $M$ is
$F$-relative injective or projective if and only if it is injective
or projective as an $A$-module. Thus $F$ is Maschke if and only if
every $A$-module is injective, and $F$ is dual Maschke if and only
if every $A$-module is projective. It is well-known that both
conditions are equivalent to $A$ being semisimple, see e.g.
\cite[Th. 5.3.7]{Cohn}. The classical Maschke Theorem can therefore
be restated as follows in our terminology: for a finite group $G$,
the restriction of scalars functor $F:\ {\cal M}_{kG}\to {\cal M}_{k}$
is a Maschke functor if and only if the order of $G$ does not divide
the characteristic of $k$. We will come back to this in \prref{4.3}.

\end{examples}

\begin{proposition}\label{pop}\prlabel{3.3}
Any $H$-separable functor $F:\ \Cc \to \Dd$ is at the same time
an $H$-Maschke and a dual $H$-Maschke functor.
\end{proposition}

\begin{proof}
We will prove first that $F$ is a $H$-Maschke functor, and leave
the proof of the second statement to the reader. Take an
object $M\in \Cc$, and let $i$ and $f$ be as in \deref{3.1}.
Then define
$$g=H(f)\circ \Pp_{C',C}(p)$$
where $p$ is a left inverse of $F(i):\ F(C)\to F(C')$ in $\Dd$.
Using \eqref{2.2.1}, we obtain
\begin{eqnarray*}
g\circ H(i)&=&H(f)\circ \Pp_{C',C}(p)\circ H(i)=
\Pp_{C,M}(F(f)\circ p\circ F(i))\\
&=& \Pp_{C,M}(F(f))=H(f)
\end{eqnarray*}
as needed.
\end{proof}

Our next result is a Rafael-type Theorem for Maschke and dual Maschke
functors.

\begin{theorem}\thlabel{rafmas}\thlabel{3.4}
Assume that the functor $F:\ \Cc\to \Dd$ has a right adjoint
$G:\ \Dd\to \Cc$.\\
1) $M\in \Cc$ is $F$-relative $H$-injective if and only if
$H(\eta_M):\ H(M) \to HGF(M)$ has a
left inverse in $\Ee$. In particular, $F$ is an $H$-Maschke functor
if and only if every $H(\eta_M)$
splits in $\Ee$.\\
2) $P\in \Dd$ is $G$-relative $H$-projective if and only
if $H(\varepsilon_P):\ HFG(P) \to H(P)$
has a right inverse. In particular, $G$ is a dual
$H$-Maschke functor if and only if every $H(\varepsilon_P)$
cosplits in $\Ee$.
\end{theorem}

\begin{proof}
1) Assume first that $M$ is $F$-relative $H$-injective.
Consider the unit map $\eta_M:\ M\to GF(M)$ in $\Cc$.
From \eqref{2.5.1}, we know that $F(\eta_M)$ has a left
inverse in $\Dd$, so there exists a map
$\nu_M:\ HGF(M)\to H(M)$ in $\Ee$ making the diagram
$$\begin{diagram}
 H(M)       &  \rTo^{H(\eta_M)}  &   HGF(M)  \\
\dTo^{I_{H(M)}} &  \SW_{\nu_M} \\
H(M)        &          &
\end{diagram}$$
commutative. This means that $\nu_M$ is a left inverse of $H(\eta_M)$.\\
Conversely, assume that $H(\eta_M)$ has a left inverse $\nu_M$,
and consider $i:\ C\to C'$, $f:\ C\to M$, with $p:\ F(C')\to F(C)$
a left inverse of $F(i)$. Then take
$$g=\nu_M\circ HGF(f)\circ HG(p)\circ H(\eta_{C'}):\
H(C')\to H(M)$$
$\eta$ is a natural transformation, hence the diagrams
$$
\begin{diagram}
 C  & \rTo^{i} & C'\\
\dTo^{\eta_C}& &\dTo_{\eta_{C'}}\\
GF(C)&\rTo^{GF(i)}&GF(C')
\end{diagram} \quad \quad \quad \quad \quad \quad \quad
\begin{diagram}
 C  & \rTo^{f} & M\\
\dTo^{\eta_C}& &\dTo_{\eta_{M}}\\
GF(C)&\rTo^{GF(f)}&GF(M)
\end{diagram}
$$
commute. Using this, we compute
\begin{eqnarray*}
g\circ H(i)&=&
\nu_M\circ HGF(f)\circ HG(p)\circ H(\eta_{C'})\circ H(i)\\
&=& \nu_M\circ HGF(f)\circ HG(p)\circ HGF(i)\circ H(\eta_C)\\
&=& \nu_M\circ HGF(f)\circ HG(p\circ F(i))\circ H(\eta_C)\\
&=&\nu_M\circ HGF(f)\circ H(\eta_C)\\
&=& \nu_M\circ H(\eta_M)\circ H(f)=H(f)
\end{eqnarray*}
and this proves that $M$ is $F$-relative $H$-injective.
The proof of 2) is left to the reader.
\end{proof}

\begin{remark}\relabel{3.5}\rm
Let us compare the Rafael Theorems for Maschke functors and
separable functors. A functor $F$ with a right adjoint $G$ is
Maschke if and only if every unit morphism $\eta_M$ has a
left inverse $\nu_M$. $F$ is separable if, moreover, $\nu$ is
natural in $M$. We have similar interpretations for $H$-Maschke
and $H$-separable functors. Also note that the two Rafael Theorems
\ref{th:2.7} and \ref{th:3.4} imply one statement of \prref{3.3}.
In the next example, we will see that a Maschke functor is not
necessarily separable.
\end{remark}

\begin{example}\exlabel{3.6}
\rm
Let $K\subset L$ be a finite purely inseparable field extension.
The restrictions of scalars functor $F :\ {\cal M}_L \to {\cal M}_K$
is a Maschke and a dual Maschke functor, since every $L$-vector
space is an injective and projective object of ${\cal M}_L$;
$F$ is not a separable functor, since $L/K$ is
not separable.
\end{example}

Now let $F:\ \Cc\to \Dd$ and $H:\ \Cc\to \Ee$ be functors
between abelian categories. We say that $F$ is {\em $H$-semisimple}
if the following assertion holds:
if we have an exact sequence
$$0\to C'\to C\to C''\to 0$$
in $\Cc$ such that
$$0\to F(C')\to F(C)\to F(C'')\to 0$$
is split exact in $\Dd$, then
$$0\to H(C')\to H(C)\to H(C'')\to 0$$
is split exact in $\Ee$. $F$ is called {\em semisimple} if it is
$1_{\Cc}$-semisimple. Our terminology is inspired by the fact
that an algebra
$A$ over a field
$k$ is semisimple if and only if the restriction of scalars functor
$\Mm_A\to \Mm_k$ is semisimple. It is now easy to prove the following
result.

\begin{proposition}\prlabel{3.7}
Let $F:\ \Cc\to \Dd$ and $H:\ \Cc\to \Ee$ be functors
between abelian categories.
1) If $F$ is $H$-Maschke, then $F$ is also $H$-semisimple; if $F$
reflects monomorphisms, then the converse is also true.\\
2) If $F$ is dual $H$-Maschke, then $F$ is also $H$-semisimple; if $F$
reflects monomorphisms, then the converse is also true.
\end{proposition}

\section{Examples and applications}\selabel{4}
\subsubsection*{Extension and restriction of scalars}
We consider ring morphisms $Q\to R\to S$ and $T\to S$. Associated
to these morphisms are the restriction of scalars functors
$$
\matrix{
\Mm_S\stackrel{G}{\longrightarrow} \Mm_R &\Mm_S\stackrel{G_1}
{\longrightarrow}
\Mm_T& \Mm_R\stackrel{G_2}{\longrightarrow} \Mm_Q\cr
{}_S\Mm \stackrel{G'}{\longrightarrow} {}_R\Mm & {}_S\Mm
\stackrel{G'_1}{\longrightarrow} {}_T\Mm &
{}_R\Mm \stackrel{G'_2}{\longrightarrow} {}_Q\Mm\cr}
$$
and their left adjoints, the induction functors
$F= -\ot_R S$, $F_1= -\ot_T S$, $F_2= -\ot_Q R$,
$F'= S\ot_R -$, $F'_1= S\ot_T -$, $F'_2= R\ot_Q -$.
With this notation, we have:

\begin{proposition}\prlabel{4.1}
The following assertions are equivalent:
\begin{itemize}
\item $G:\ \Mm_S\to \Mm_R$ is $G_1$-separable;\\
\item $G':\ {}_S\Mm\to {}_R\Mm$ is $G'_1$-separable;\\
\item there exists an element $e=\sum e^1\ot_R e^2\in S\ot_R S$ such that
\end{itemize}
\begin{eqnarray}
&&\sum te^1\ot_R e^2=\sum e^1\ot_R e^2t,~~~{\rm for~all}~t\in T
\eqlabel{4.1.1}\\
&&\sum e^1e^2=1\eqlabel{4.1.2}
\end{eqnarray}
\end{proposition}

\begin{proof}
Basically, this follows from the fact that $\Nat(G_1,G_1FG)$ is
in bijective correspondence with the set of $e$ satisfying
\eqref{4.1.1}: for a natural transformation $\zeta:\ G_1\to G_1FG$,
the map $\zeta_S:\ S\to S\ot_R S$ is right $T$-linear. For any
$a\in S$, consider $f_a:\ S\to S$, $f_a(s)=as$. Then $f_a\in \Mm_S$,
and the naturality of $\zeta$ implies that
$$(f_a\ot I_S)(\zeta_S(s))=\zeta_S(f_a(s))$$
Let $s=1$ and $\zeta_S(1)=\sum e^1\ot_R e^2$. Then \eqref{4.1.1} follows.
Conversely, given $e$ satisfying \eqref{4.1.1}, we construct a
natural transformation $\zeta$ as follows:
$$\zeta_M:\ M\to M\ot_RS;~~\zeta_M(m)=\sum me^1\ot_R e^2$$
It follows from \eqref{4.1.1} that $\zeta_M$ is right $T$-linear, and
we leave it to the reader to show that $\zeta$ is natural.\\
If $e$ satisfies \eqref{4.1.2}, then for all $M\in \Mm_R$ and
$m\in M$:
$$(G_1(\varepsilon_M)\circ\zeta_M)(m)=
\varepsilon_M(\sum me^1\ot_R e^2)=\sum me^1 e^2=m$$
and it follows from \thref{2.7} that $G$ is $G_1$-separable. The converse,
and the equivalence between the second and third
assertion is done in a similar way.
\end{proof}

Let us explain what this means. It is well-known \cite{NastasescuVV89},
and actually a special case of \prref{4.1}, that
$\Mm_S\to \Mm_R$ is a separable
functor if and only if $R\to S$ is separable
in the sense of \cite{DI71}, which means that there exists $e\in S\ot_RS$
satisfying \eqref{4.1.2}, and \eqref{4.1.1} also, but for all $t\in S$.
In this situation, an exact sequence in $\Mm_S$ that splits in $\Mm_R$
also splits in $\Mm_S$.
In \prref{4.1}, we have $e\in S\ot_RS$
satisfying \eqref{4.1.2}, and \eqref{4.1.1}, but only for $t$ in a subring
$T$ of $S$. We then have the weaker conclusion that an
exact sequence in $\Mm_S$ that splits in $\Mm_R$
also splits in $\Mm_T$.\\
On the other hand, a nice ring-theoretical
problem arises from the concept of (dual) $H$-Maschke
functor:

{\em Let $T\to S$ be a ring morphism. When is any right $S$-module
is projective (injective) as a right $T$-module?}

Using \prref{4.1}, we obtain a suficient condition is obtain:

\begin{corollary}
Let $T\to S$ be a morphism of $k$-algebras over a field $k$
such that $S$ is projective as a right $T$-module. Assume that
there exists $e = \sum e^1 \ot e^2 \in (S\ot S)^T$ such that
$\sum e^1 e^2 =1_S$. Then any  right $S$-module
is projective as a right $T$-module.
\end{corollary}

\begin{proof}
We take $R=k$ is the \prref{4.1}. If such an $e$ exists, then
the forgetful functor $G: {\cal M}_S \to {\cal M}_k$ is
$G_1$-separable, where $G_1 :{\cal M}_S \to {\cal M}_T$ is the
restriction of scalars functor. Hence, $G$ is a dual $G_1$-Maschke
functor. Let $M$ be a right $S$-module; as $k$ is a field, the
right $S$-module structure on $M$, $\nu_M :M\ot S \to M$, has a
section in ${\cal M}_k$. Thus, there exists
$f: M\to M\ot S \cong S^{(M)}$ a right $T$-module map such that
$\nu_M \circ f ={\rm Id}_M$, i.e. $M$ is a direct summand of
$S^{(M)}$ as a right $T$-submodule. As $S$ is projective in
${\cal M}_T$ we obtain that $M$ is projective as a right
$T$-module.
\end{proof}

We have a similar result for split extensions:

\begin{proposition}\prlabel{4.2}
The induction functor $F= -\ot_R S$ is $G_2$-separable if and only if the
ring morphism $R\to S$ is split as a map of $(R,Q)$-bimodules.
\end{proposition}

\begin{proof}
Assume that $F$ is $G_2$-separable. According to \thref{2.7}, there
exists a natural transformation $\nu:\ G_2GF\to G_2$ such that
$$\nu_M\circ G_2(\eta_M)=I_{G_2(M)}$$
for all $M\in \Mm_R$. This means that $\nu_R:\ S\to R$ splits
$R\to S$ as a map of right $Q$-modules. From the naturality of $\nu$,
we can deduce that $\nu_R$ is also left $R$-linear, and it follows
that $R\to S$ is split as a map of $(R,Q)$-bimodules. The converse is
left to the reader.
\end{proof}

We now assume that $(F,G)$ is a Frobenius pair of functors; this means
that the ring extension $R\to S$ is Frobenius, and it is equivalent
to the existence of a Frobenius system (cf. e.g.
\cite{BrzezinskiCMZ00} or \cite{Kadison99b}).
A Frobenius system consists of a pair $(\ol{\mu},f)$,
where $\ol{\mu}:\ S\to R$ is an $R$-bimodule map,
$f=\sum f^1\ot_R f^2\in S\ot_R S$ is a Casimir element
i.e.
$$\sum sf^1\ot_R f^2=\sum f^1\ot_R f^2s$$
for all $s\in S$ and
\begin{equation}\eqlabel{4.2b.1}
\sum \ol{\mu}(f^1)f^2=\sum f^1\ol{\mu}(f^2)=1
\end{equation}
Our next two results can be deduced from \prref{2.9}, but it is easier
to give a direct proof.

\begin{proposition}\prlabel{4.2a}
We keep the notation from above, assuming that the ring
extension $R\to S$ is Frobenius, with Frobenius system $(\ol{\mu},f)$.
Then $G$ is $G_1$-separable if and only if there exists an
$(R,T)$-bimodule map $\alpha:\ S\to S$ such that
$\sum f^1\alpha(f^2)=1$.
\end{proposition}

\begin{proof}
Assume that $G$ is $G_1$-separable, and take
$e=\sum e^1\ot_R e^2 \in S\ot_R S$
as in \prref{4.1}. We define $\alpha:\ S\to S$ by
$$\alpha(s)=\sum \ol{\mu}(se^1)e^2=\sum \ol{\mu}(e^1)e^2s$$
Using the fact that $\ol{\mu}$ is left $R$-linear, we easily prove
that $\alpha$ is left $R$-linear. For all $s\in S$ and $t\in T$, we have
$$\alpha(st)=\sum \ol{\mu}(ste^1)e^2=\sum \ol{\mu}(se^1)e^2t=\alpha(s)t$$
so $\alpha$ is right $T$-linear. Finally
$$\sum f^1\alpha(f^2)=\sum f^1\ol{\mu}(f^2e^1)e^2=
\sum e^1f^1\ol{\mu}(f^2)e^2=\sum e^1e^2=1$$
Conversely, suppose that we have an $(R,T)$-bimodule map
$\alpha : S\to S$
such that $\sum f^1\alpha(f^2)=1$. We then take
$$e=\sum e^1\ot_R e^2=\sum f^1\ot_R \alpha(f^2)\in S\ot_R S$$
and compute that $\sum e^1e^2=1$ and
$$\sum tf^1\ot_R \alpha(f^2)=
\sum f^1\ot_R \alpha(f^2t)=\sum f^1\ot_R \alpha(f^2)t$$
and it follows from \prref{4.1} that $G$ is $G_1$-separable.
\end{proof}

\begin{proposition}\prlabel{4.2b}
We keep the notation from above, assuming that the ring
extension $R\to S$ is Frobenius, with Frobenius system $(\ol{\mu},f)$.
Then the following assertions are equivalent:
\begin{itemize}
\item $F= -\ot_R S$ is $G_2$-separable;
\item $F= -\ot_R S$ is $G'_2$-separable;
\item there exists $x\in C_Q(S)$
such that $\ol{\mu}(x)=1$.
\end{itemize}
\end{proposition}

\begin{proof}
First assume that $F$ is $G_2$-separable. From \prref{4.2}, we know
that there exists an $(R,Q)$-bimodule map $\ol{\nu}:\ S\to R$
such that $\ol{\nu}(1_S)=1_R$. Take $x=\sum f^1\ol{\nu}(f^2)$.
Then for all $q\in Q$, we have
$$qx=\sum qf^1 \ol{\nu}(f^2)=\sum f^1 \ol{\nu}(f^2q)
=\sum f^1 \ol{\nu}(f^2)q=xq$$
and
\begin{eqnarray*}
\ol{\mu}(x)&=& \sum \ol{\mu}\bigl(f^1\ol{\nu}(f^2)\bigr)=
\sum \ol{\mu}(f^1)\ol{\nu}(f^2)\\
&=&\sum \ol{\nu}\bigl(\ol{\mu}(f^1)f^2\bigr)=\ol{\nu}(1_S)=1_R
\end{eqnarray*}
Conversely, given $x\in C_Q(S)$
such that $\ol{\mu}(x)=1$, we define $\ol{\nu}:\ S\to R$ by
$\ol{\nu}(s)=\ol{\mu}(sx)$. Then $\ol{\nu}(1)=\ol{\mu}(x)=1$,
and, for all $r\in R$, $s\in S$ and $q\in Q$, we have
$$\ol{\nu}(rs)=\ol{\mu}(rsx)=r\ol{\mu}(sx)=r\ol{\nu}(s)$$
$$\ol{\nu}(sq)=\ol{\mu}(sqx)=\ol{\mu}(sxq)=\ol{\mu}(sx)q
=\ol{\nu}(s)q$$
and $\ol{\nu}$ is an $(R,Q)$-bimodule map, as needed.\\
The equivalence between the second and the third assertion can be
shown in a similar way.
\end{proof}

\begin{remark}\relabel{4.2c}\rm
It follows from \prref{4.1} that the
$G_1$-separability of the restriction of scalars functor $G$ is
left-right symmetric. It is remarkable that a similar property does
not hold for the $G_2$-separability of the induction functor
(see \prref{4.2}), unless we know that $S/R$ is Frobenius
(see \prref{4.2b}).
\end{remark}

\subsubsection*{Hopf algebras}
A separable functor is always Maschke and dual Maschke, but the
converse is in general not true, see \exref{3.6}. However, there
are some particular situations where the converse property holds.\\
A classical result of Sweedler (\cite{Sweedler69}) states
that a Hopf algebra over a field is semisimple if and only if there
exists a (left or right) integral $t\in H$ such that
$\varepsilon(t)=1$. The generalization to Hopf algebras over a commutative
ring $k$ is the following: a Hopf algebra is separable if and only if
there exists an integral $t$ with $\varepsilon(t)=1$. Here the remarkable
thing is that, over a field $k$, a separable algebra is semisimple,
but not conversely: it suffices to look at a purely inseparable field
extension. We can now explain this apparent contradiction. First
observe the following.\\
An algebra $A$ over a field $k$ is semisimple if and only if the
restriction of scalars functor $\Mm_A\to \Mm_k$ is a Maschke functor,
if and only if it is a dual Maschke functor. This is a restatement of
the classical result \cite[Th. 5.3.7]{Cohn}.\\
An algebra $A$ over a commutative ring $k$ is separable if and only if
$\Mm_A\to \Mm_k$ is a separable functor (see \prref{4.1} with $T=S$ or
\cite{NastasescuVV89}).\\
With these observations in mind, we restate and prove Sweedler's results
in the following fashion.

\begin{proposition}\prlabel{4.3}
Let $H$ be a Hopf algebra over a commutative ring $k$, and
$G:\ \Mm_H\to \Mm_k$ the restriction of scalars functor. Then
the following
assertions are equivalent:\\
1) $G$ is a
dual Maschke functor;\\
2) $G$ is a
Maschke functor;\\
3) $G$ is a
semisimple functor;\\
4) there exists a right integral $t\in H$ with $\varepsilon(t)=1$;\\
5) $G$ is a
separable functor.
\end{proposition}

\begin{proof}
The equivalence of 1), 2) and 3) follows immediately from \prref{3.7},
since $G$ reflects monomorphisms and epimorphisms.\\
$1)~\Rightarrow~ 4)$. $k\in \Mm_H$, with the trivial action:
$x\cdot h=\varepsilon(h)x$. Then $\varepsilon:\ H\to k$ in $\Mm_H$
is such that $G(\varepsilon)$ is a cosplit epimorphism. So we have
a map $\tau\in \Mm_H$ making the following diagram commutative
in $\Mm_H$:
$$\begin{diagram}
H&\rTo^{\varepsilon}& k\\
\dTo_{I_H}&\SW^{\tau}&\\
H&&\end{diagram}$$
$t=\tau(1)$ is then the required integral.\\
$4)~\Rightarrow~ 5)$. Let $t$ be a right integral, with $\varepsilon(t)=
1$. $S(t_{(1)})\ot t_{(2)}$ is the required separability idempotent.\\
$5)~\Rightarrow~ 1)$ follows from \prref{3.3}.
\end{proof}

The dual version of this result is the following; we leave the proof
to the reader.

\begin{proposition}\prlabel{4.4}
Let $H$ be a flat Hopf algebra over a commutative ring $k$, and
$F:\ \Mm^H\to \Mm_k$ the forgetful functor. Then
the following
assertions are equivalent:\\
1) $F$ is a
Maschke functor;\\
2) $F$ is a
dual Maschke functor;\\
3) $F$ is a
semisimple functor;\\
4) there exists a right integral $\varphi\in H^*$ with
$\varphi(1)=1$;\\
5) $F$ is a
separable functor.
\end{proposition}

\subsubsection*{The structure of injective objects in the
category of entwined modules}
As an application of \thref{3.4}, we give the structure of
injective objects in the category of entwined modules.
For other results about injective objects in
the category of graded modules and the category
of modules graded by a $G$-set (which are special cases of
entwined modules), we refer
to  \cite[Th. 2.1, Cor. 2.2]{Nastasescu} and
\cite[Cor. 3.3]{NastasescuRV90}; other results for Doi-Koppinen
Hopf modules in general can be found in \cite[Cor. 2.9]{CaenepeelMZ97a} and
\cite[Th. 4.3]{DNT}.\\
Let $(A,C,\psi)$ be an entwining structure of a commutative ring $k$.
We will assume that $C$ is flat as a $k$-module, to ensure that the
category of $C$-comodules and the category of entwined modules is
abelian. We will use the following notation for functors forgetting
actions and coactions:
\begin{equation}\eqlabel{4.5.1}
\begin{diagram}
\Mm(\psi)_A^C&\rTo^F&\Mm_A\\
\dTo_H&&\dTo^{H_1}\\
\Mm^C&\rTo^{F_1}&\Mm_k
\end{diagram}
\end{equation}
$F$ has a right adjoint $G=\bullet \ot C$, and $H_1$ has a right
adjoint $K_1$ given by
$$K_1(V)=\Hom(A,V),~~{\rm with}~~(f\cdot a)(b)=f(ab)$$
for any $f:\ A\to V$. Thus we have also an adjoint pair $(H_1F, GK_1)$
and the unit and counit of this adjoint pair are
$$\eta_M:\ M\to \Hom(A,M)\ot C,~~\eta_M(m)=m_{[0]}\bullet\ot m_{[1]}$$
$$\varepsilon_V:\ \Hom(A,V)\ot C\to V,~~\varepsilon_V(f\ot c)=
f(1_A)\varepsilon_C(c)$$
For any $m\in M\in \Mm_A$, we write $m\bullet$ for the map $A\to M$,
sending $a$ to $ma$.

\begin{corollary}\colabel{4.5}
Let $(A,C,\psi)$ be an entwining structure over a field $k$.
$Q$ is an injective object in $\Mm(\psi)_A^C$
if and only if there exists a vector space $V$ such that
$Q$ is isomorphic to a direct summand of $\Hom(A, V)\ot C$.
\end{corollary}

\begin{proof}
As $k$ is a field, then the category $\Mm(\psi)_A^C$ is
Grothendieck (see \cite{book}).
The forgetful functor $H_1F$ is exact and
${\cal M}_k$ has enough injectives, so the right
adjoint $GK_1$ preserves injectives. Thus
$\Hom(A, V)\ot C$ is an injective object of
${\cal M}(\psi)_A^C$, and so are its direct summands.\\
Conversely, assume that $Q$ is an injective object of
${\cal M}(\psi)_A^C$. As $k$ is a field, $Q$ is $F$-relative
injective, and it follows from \thref{3.4} that the unit
$\eta_Q :\ Q \to \Hom (A, Q)\ot C$ has a retraction in
the Grothendieck category ${\cal M}(H)_A^C$, and this means that
$Q$ is isomorphic to a direct summand of
$\Hom (A, Q)\ot C$.
\end{proof}

Let us present some examples, where the entwining structure comes
from a Doi-Koppinen datum $(H, A, C)$.

\begin{examples}\exlabel{4.6}
\rm
1. Let $(H, A, C)= (k, A, k)$; then
${\cal M}(k)_A^k= {\cal M}_A$, the category of right
$A$-modules. From \coref{4.5}, we recover the
well-known result stating that a right $A$-module $Q$ is
injective if and only if there exists a vector space
$V$ such that $Q$ is a direct summand of the right
$A$-module $\Hom (A, V)$.

2. Now let $(H, A, C)= (k, k, C)$; then
${\cal M}(k)_k^C= {\cal M}^C$, the category of
right $C$-comodules, \coref{4.5} tells that the injective
right $C$-comodule are the direct summands of $C^{(I)}$,
with $I$ an index set.

3. \coref{4.5} can be used to describe injective
modules graded by $G$-sets: let $G$ be a group, $X$ is a right $G$-set,
$A$ a $G$-graded $k$-algebra,
and consider the Doi-Koppinen datum
$(H, A, C) = (kG, A, kX)$. The corresponding Doi-Koppinen Hopf
modules are then exactly the $A$-modules graded by $X$, as
introduced in \cite{NastasescuRV90}, and it follows that
the injective objects in the category of $A$-modules graded by $X$
are the direct summands of $A$-modules graded by $X$ of the form
$\Hom (A, V)^{(X)} = \oplus_{x\in X}\  \Hom (A, V)_x$,
with $\Hom (A, V)_x= \Hom (A, V)$ for all $x\in X$.
\end{examples}

\begin{remark}\relabel{4.7}
\rm
If the forgetful functor
$F:\ {\cal M}(\psi)_A^C \to {\cal M}_k$
has a left adjoint, then we can also describe
the projective objects of ${\cal M}(\psi)_A^C$.
Unfortunately, in general, $F$ has not a left
adjoint: for instance, the forgetful functor
$F :\ {\cal M}^C \to {\cal M}_k$ has
a left adjoint if and only if $C$ is finite dimensional
over the field $k$ (this result is a special case of
\cite[Proposition 1.10]{Takeuchi77}); in this
case ${\cal M}^C \cong {}_{C^*}{\cal M}$, the
category of modules over $C^*$.
\end{remark}

\subsubsection*{Entwined modules and separability}
We keep the notation \eqref{4.5.1}. Let us examine when $F$ is
$H$-separable. In order to apply \thref{2.7}, we need to
examine $\Nat(HGF,H)$. In \cite[Proposition 4.1]{BrzezinskiCMZ00},
$\Nat(GF,1)$ has been computed, and an adaption of the
arguments leads to a description of $\Nat(HGF,H)$. We present
a brief sketch: consider a natural transformation
$\nu:\ HGF\to H$. $A\ot C=G(A)\in \Mm(\psi)_A^C$, so we can
consider the map
$$\nu_{A\ot C}:\ HGF(A\ot C)=A\ot C\ot C\to H(A\ot C)=A\ot C$$
in $\Mm^C$. Now we define $\theta:\ C\ot C\to A$ by
$$\theta(c\ot d)=(I_A\ot \varepsilon)(\nu_{A\ot C}(1\ot c\ot d))$$
Using the naturality of $\nu$, we can prove that $\theta$
satisfies the relation
\begin{equation}\eqlabel{4.8.1}
\theta(c\ot d_{(1)})\ot d_{(2)}=\theta(c_{(2)}\ot d)_{\psi}
\ot c_{(1)}^{\psi}
\end{equation}
for all $c$, $d\in C$.
Conversely, given a map $\theta:\ C\ot C\to A$ satisfying \eqref{4.8.1},
we can define a natural transformation $\nu:\ HGF\to H$ as follows:
let
$$\nu_M:\ M\ot C\to M~~;~~\nu_M(m\ot c)=m_{[0]}\theta(m_{[1]}\ot c)$$
for all $M\in \Mm(\psi)_A^C$.
It is clear that $\nu_M\circ H(\eta_M)=I_{H(M)}$,
for all $M\in \Mm(\psi)_A^C$
if and only if $\theta(\Delta_C(c))=\varepsilon_C(c)1_A$,
for all $c\in C$. If such a map $\theta$ exists, then
$\nu_M$ is a retraction in $\Mm^C$ of the $C$-coaction
$\eta_M = \rho_M : M\to M\ot C$. Thus any $M\in \Mm(\psi)_A^C$
is relative injective as a right $C$-comodule.
We summarize our result in the following Proposition which is
an equivalent version for entwining modules of
\cite[Theorem 2.6]{MeniniMilitaru00a}.

\begin{proposition}\prlabel{4.8}
Let $(A,C,\psi)$ be an entwining structure, and consider the forgetful
functors $F:\ \Mm(\psi)_A^C\to \Mm_A$ and $H:\ \Mm(\psi)_A^C\to \Mm^C$.
Then $F$ is $H$-separable if and only if there exists a map
$\theta:\ C\ot C\to A$ such that
\begin{equation}\eqlabel{4.8.2}
\theta(c\ot d_{(1)})\ot d_{(2)}=\theta(c_{(2)}\ot d)_{\psi}
\ot c_{(1)}^{\psi}  \qquad {\rm and} \qquad
\theta\circ\Delta_C = \eta_A\circ \varepsilon_C
\end{equation}
for all $c$, $d\in C$. In this case any $M\in \Mm(\psi)_A^C$
is relative injective as a right $C$-comodule.
\end{proposition}

In a similar way, we can investigate when the functor $H$ is
$F$-separable. We then obtain the following:

\begin{proposition}\prlabel{4.9}
Let $(A,C,\psi)$ be an entwining structure, and consider the forgetful
functors $F:\ \Mm(\psi)_A^C\to \Mm_A$ and $H:\ \Mm(\psi)_A^C\to \Mm^C$.
Then $H$ is $F$-separable if and only if there exists a map
$$e:\ C\to A\ot A,~~~e(c)=\sum e^1(c)\ot e^2(c)$$
such that
\begin{equation}\eqlabel{4.9.1}
\sum e^1(c)\ot e^2(c)a=
\sum a_{\psi}e^1(c^{\psi})\ot e^2(c^{\psi}) \quad {\rm and} \quad
m_A\circ e = \eta_A\circ \varepsilon_C
\end{equation}
for all $c\in C$, $a\in A$. In this case any $M\in \Mm(\psi)_A^C$
is relative projective as a right $A$-module.
\end{proposition}

Recall from \cite[Theorem 3.4]{BrzezinskiCMZ00} that the pair $(F,G)$ is
Frobenius if and only if there exists a $k$-linear map $\theta:\ C\ot C\to A$
and $z=\sum_l a_l\ot c_l\in A\ot C$ such that
\begin{eqnarray*}
&&\theta(c\ot d)a= a_{\psi\Psi}\theta(c^{\Psi}\ot d^{\psi})\\
&&\theta(c\ot d_{(1)})\ot d_{(2)}=
\theta(c_{(2)}\ot d)_{\psi}\ot c_{(1)}^{\psi}\\
&&az=za\\
&&\eta_A(\varepsilon_C(d))=\sum_l a_l\theta(c_l\ot d)=
\sum_l a_{l\psi}\theta(d^{\psi}\ot c_l)
\end{eqnarray*}
for all $a\in A$ and $c\in C$. We will call $(\theta,z)$ a Frobenius
system for the adjunction $(F,G)$. We give the unit and counit of
the adjunctions $(F,G)$ and $(G,F)$:
\begin{eqnarray*}
\eta:\ 1\to GF~~&~~\eta_M:\ M\to M\ot C~~&~~\eta_M(m)=m_{[0]}\ot m_{[1]}\\
\varepsilon:\ FG\to 1~~&~~\varepsilon_N:\ N\ot C\to N~~&~~
\varepsilon_N(n\ot c)=\varepsilon_C(c)n\\
\nu:\ GF\to 1~~&~~\nu_M:\ M\ot C\to M~~&~~\nu_M(m\ot c)=
m_{[0]}\theta(m_{[1]}\ot c)\\
\zeta:\ 1\to FG~~&~~\zeta_N(n)=\sum_l na_l\ot c_l
\end{eqnarray*}
Now assume that $(F,G)$ is a Frobenius pair, and that we know a Frobenius
system $(\theta,z)$. Using \prref{2.9}, we can decide when $F$ is
separable or $H$-separable.

\begin{lemma}\lelabel{4.9a}
With notation as above,
$$\Nat(HG,HG)\cong \Hom(C, A)$$
and
$$\Nat(G,G)\cong \{\beta\in \Hom(C, A)~|~\beta(c)a=a_{\psi}\beta(c^{\psi})~
{\rm for~all}~a\in A,~c\in C\}$$
\end{lemma}

\begin{proof}
Consider a natural transformation $\alpha:\ HG\to HG$. Then the
map $\alpha_A:\ A\ot C\to A\ot C$ is right $C$-colinear, and, using the
naturality of $\alpha$, we find that $\alpha_A$ is also left $A$-linear.
Now consider the map $\beta:\ C\ot A$ defined by
$$\beta(c)=(I_A\ot \varepsilon_C)(\alpha_A(1_A\ot c))$$
Conversely given $\beta:\ C\to A$, we define a natural transformation
$\alpha:\ HG\to HG$ by putting
$$\alpha_N:\ N\ot C\to N\ot C,~~\alpha_N(n\ot c)=n\beta(c_{(1)})
\ot c_{(2)}$$
for every $N\in \Mm_A$. It is obvious that $\alpha_N$ is right $C$-colinear;
let us check that $\alpha$ is natural. For all $f:\ N\to N'$ in
$\Mm_A$, we have
\begin{eqnarray*}
&&\hspace*{-2cm}\alpha_{N'}(f(n)\ot c)=f(n)\beta(c_{(1)})\ot c_{(2)}\\
&=&f(n\beta(c_{(1)}))\ot c_{(2)}=(f\ot I_C)(\alpha_N(n\ot c))
\end{eqnarray*}
If $\alpha:\ G\to G$ is a natural transformation, then the map
$\alpha_A$ is also right $A$-linear, and it follows easily that
$\beta$ defined as above satisfies the centralizing condition
\begin{equation}\eqlabel{4.9a.1}
\beta(c)a=a_{\psi}\beta(c^{\psi})
\end{equation}
If $\beta:\ C\to A$ satisfies \eqref{4.9a.1}, then we define
$\alpha:\ G\to G$ by the same formula as above, and the second statement
of the Lemma follows if we can prove that $\alpha_N$ is right $A$-linear.
This goes as follows:
\begin{eqnarray*}
&&\alpha_N((n\ot c)a)=\alpha_N(na_{\psi}\ot c^{\psi})\\
&=& na_{\psi}\beta((c^{\psi})_{(1)})\ot  c^{\psi})_{(2)}\\
{\rm \eqref{1.3}}~~~~&=& na_{\psi\Psi}(\beta(c_{(1)}^{\Psi}))\ot
c_{(2)}^{\psi})\\
 {\rm \eqref{4.9a.1}}~~~~&=& n\beta(c_{(1)})a_{\psi}\ot c_{(2)}^{\psi})\\
&=& (n\beta(c_{(1)})\ot c_{(2)})a=\alpha_N(n\ot c)a
\end{eqnarray*}
\end{proof}

\begin{theorem}\thlabel{4.9b}
Consider the forgetful functors $F:\ \Mm(\psi)_A^C\ot \Mm_A$
and $H:\ \Mm(\psi)_A^C\ot \Mm^C$, and assume that the functor $F$ and
its adjoint form a Frobenius pair, with Frobenius system $(\theta,z)$. 
Then $F$ is $H$-separable if and only if there exists a a map
$\beta:\ C\ot A$ such that
\begin{equation}\eqlabel{4.9b.1}
\beta(c_{(3)})_{\psi_1\psi_2}\theta(c_{(2)}^{\psi_1}\ot c_{(4)})_{\psi_3}
\ot c_{(1)}^{\psi_2\psi_3}=1\ot c
\end{equation}
for all $c\in C$. $F$ is separable if and only if there exists a
$\beta:\ C\ot A$ satisfying \eqref{4.9b.1} and \eqref{4.9a.1}.
\end{theorem}

\begin{proof}
This follows immediately from \prref{2.9}: if we apply \eqref{2.9.1}
to $A\in \Mm_A$, then we find \eqref{4.9b.1}. Conversely, if
$\beta$ satisfies \eqref{4.9b.1}, then we easily compute that the
corresponding natural transformation $\alpha$ satisfies \eqref{2.9.1}.
\end{proof}

Let us now discuss the dual version of \thref{4.9b}. Let $K$ be
the left adjoint of the forgetful functor $H:\ \Mm(\psi)_A^C\to
\Mm^(\psi)$. In \cite[Proposition 4.4]{BrzezinskiCMZ00}, it is shown
that $(H,K)$ is a Frobenius pair if and only if there exists a
Frobenius system $(\vartheta, e)$, consisting of maps
$\vartheta\in (C\ot A)^*$ and $e:\ C\to A\ot A$ such that
\begin{eqnarray*}
&&\vartheta(c_{(1)}\ot a_{\psi})c_{(2)}^{\psi}=
\vartheta(c_{(2)}\ot a)c_{(1)}\\
&&e^1(c_{(1)})\ot e^2(c_{(1)})\ot c_{(2)}=
e^1(c_{(2)})_{\psi}\ot e^2(c_{(2)})_{\Psi}\ot c_{(1)}^{\psi\Psi}\\
&&e^1(c)\ot e^2(c)a= a_{\psi}e^1(c^{\psi})\ot e^2(c^{\psi})\\
&&\varepsilon(c) 1= \vartheta(c_{(1)}\ot e^1(c_{(2)})) e^2(c_{(2)})
= \vartheta(c_{(1)}^{\psi}\ot e^2(c_{(2)})) e^1(c_{(2)})_{\psi}
\end{eqnarray*}
for all $c\in C$ and $a\in A$. We use the notation
$$e(c)=e^1(c)\ot e^2(c)$$
with summation implicitely understood.
The unit and counit of the adjunction $(H,K)$ is then given by
\begin{eqnarray*}
\zeta:\ 1\to KH~~&~~\zeta_M:\ M\ot M\ot A~~&~~
\zeta_M(m)=m_{[0]}e^1(m_{[1]})\ot e^1(m_{[1]})\\
\nu:\ HK\to 1~~&~~\nu_N:\ N\ot A\to N~~&~~\nu_N(n\ot a)=
\vartheta(n_{[1]}\ot a)n_{[0]}
\end{eqnarray*}

\begin{lemma}\lelabel{4.9c}
We have isomorphisms
$$\Nat(FK,FK)=\Hom(C,A)$$
and
$$\Nat(K,K)=\{\beta\in \Hom(C,A)~|~
\beta(c_{(1)})\ot c_{(2)}=\beta(c_{(2)})_{\psi}\ot c_{(1)}^{\psi}~
{\rm for~all~}c\in C\}$$
\end{lemma}

\begin{theorem}\thlabel{4.9d}
Consider the forgetful functors $F:\ \Mm(\psi)_A^C\ot \Mm_A$
and $H:\ \Mm(\psi)_A^C\ot \Mm^C$, and assume that the functor $H$ and
its adjoint form a Frobenius pair, with Frobenius system $(\vartheta,e)$. 
Then $H$ is $F$-separable if and only if there exists a a map
$\beta:\ C\ot A$ such that
\begin{equation}\eqlabel{4.10b.1}
\varepsilon(c)a=a_{\psi_1\psi_2}e^1(c_{(2)}^{\psi_1})_{\psi_3}
\beta(c_{(1)}^{\psi_2\psi_3})e^2(c_{(2)}^{\psi_1})
\end{equation}
for all $c\in C$ and $a\in A$. $H$ is separable if and only if there exists a
$\beta:\ C\ot A$ satisfying \eqref{4.10b.1} and 
$$\beta(c_{(1)})\ot c_{(2)}=\beta(c_{(2)})_{\psi}\ot c_{(1)}^{\psi}$$
for all $c\in C$.
\end{theorem}

\subsubsection*{Yetter-Drinfeld modules and quantum integrals}
\prref{4.8} and \prref{4.9} can be applied in many situations:
Doi-Koppinen modules, Yetter-Drinfeld modules, relative Hopf
modules, graded modules, etc. are all special cases of the
category $\Mm(\psi)_A^C$. In this subsection we shall apply
the above results to the category ${\cal YD}^H_H$ of Yetter-Drinfeld
modules \cite{Y}.

Let $(A, C, \psi) = (L, L,\psi)$, where $L$ is a Hopf algebra
and
$$
\psi : L\ot L \to L\ot L, \quad
\psi (g\ot h) = h_{(2)} \ot S(h_{(1)}) g h_{(3)}
$$
for all $g$, $h\in L$. The resulting category of entwined modules
is just
$\Mm(\psi)_L^L = {\cal YD}^L_L$, the category of Yetter-Drinfeld
modules over $L$.

\begin{corollary}\colabel{ydint}
Let $L$ be a Hopf algebra over a commutative ring $k$ and consider
the forgetful functors
$F:\ {\cal YD}^L_L\to \Mm_L$ and $H:\ {\cal YD}^L_L\to \Mm^L$.\\
1) The following statements are equivalent:

$\bullet$ $F$ is $H$-separable;

$\bullet$ there exists a $k$-linear map $\theta :\ L\ot L \to L$
such that
$$
\theta (g \ot h_{(1)} ) \ot h_{(2)} =
\theta (g_{(2)} \ot h)_{(2)} \ot
S\Bigl(\theta (g_{(2)} \ot h)_{(1)}\Bigl) g_{(1)}
\theta (g_{(2)} \ot h)_{(3)}, \quad
\theta (h_{(1)} \ot h_{(2)}) = \varepsilon (h) 1_H
$$
for all $g$, $h\in L$;

$\bullet$ there exists a $k$-linear map $\gamma :\ L \to \End(L)$
such that
$$
\gamma (h_{(1)})(g) \ot h_{(2)} =
\gamma (h)(g_{(2)})_{(2)} \ot
S\Bigl(\gamma (h)(g_{(2)})_{(1)} \Bigl) g_{(1)}
\gamma (h)(g_{(2)})_{(3)}, \quad
\gamma (h_{(2)})(h_{(1)}) = \varepsilon (h) 1_H
$$
for all $g$, $h\in L$. In this case any $M\in {\cal YD}^L_L$
is relative injective as a right $L$-comodule.

2) The following statements are equivalent:

$\bullet$ $H$ is $F$-separable;

$\bullet$ there exists a $k$-linear map $e : L\to L \ot L$,
$e(h)=\sum e^1(h)\ot e^2(h) \in L\ot L$
such that
$$
\sum e^1(g)\ot e^2(g) h =
\sum h_{(2)} e^1 \Bigl(S(h_{(1)}) g h_{(3)} \Bigl) \ot
e^2 \Bigl(S(h_{(1)}) g h_{(3)} \Bigl), \quad
\sum e^1(g)e^2(g) = \varepsilon (g) 1_H
$$
for all $g$, $h\in L$. In this case any $M\in {\cal YD}^L_L$
is relative projective as a right $L$-module.

Furthermore, if $L$ is finitely generated and projective
over $k$, these conditions are also equivalent to

$\bullet$ There exists an element
$\sum_{i=1}^n f_i \ot h_i \in \End(L) \ot L$ such that
$$
\sum_{i=1}^n f_i(g) \ot h_i h = \sum_{i=1}^n
h_{(2)} f_i \Bigl(S(h_{(1)}) g h_{(3)} \Bigl) \ot h_i, \quad
\sum_{i=1}^n f_i(g) h_i = \varepsilon (g) 1_H
$$
for all $g$, $h\in L$.
\end{corollary}

\begin{proof}
This follows from \prref{4.8} and \prref{4.9} applied to the
above entwining structure. The equivalence between the maps
$\theta : L\ot L \to L$ and the maps $\gamma: L\to \End(L)$
is given by
the $k$-linear isomorphism given by the adjunction
$$
\Hom (L\ot L, L)\cong \Hom (L, \End(L) )
$$
Hence, for any $\theta: L\ot L \to L$ there exists a unique
$\gamma = \gamma_{\theta} : L\to \End(L)$ such that
$\theta (g\ot h) =\gamma (h) (g)$, for any $g$, $h\in L$.

In the case that $L$ is finitely generated and projective
over $k$ we use the ``Hom-tensor relations"
$$
\End(L) \ot L \cong \Hom (L, L\ot L)
$$
i.e. for any $e:\ L\to L\ot L$ there exists a unique
element $\sum_{i=1}^n f_i \ot h_i \in \End(L) \ot L$ such that
$e (g) = \sum_{i=1}^n f_i(g) \ot h_i$, for any $g\in G$.
\end{proof}

\begin{remarks}\rm
1. In \cite{MeniniMilitaru00a}, a map $\gamma :\ L \to \End(L)$
satisfying the conditions of \coref{ydint} has been called
a {\em total quantum integral}.\\
2. Assume now that $L$ is a finite dimensional Hopf algebra
over a field $k$ and let
$\sum_{i=1}^n f_i \ot h_i \in \End(L) \ot L$ be an element
as in \coref{ydint}. Let $D(L)$ be the Drinfeld double of
$L$. Then it is well know that there exists an equivalence of
categories
${\cal YD}^L_L \cong {\cal M}_D(L)$ and the above functor $F$
is just the restriction of scalars. From ring theoretical point
of view the ring extension $D(L)/L$ has a remarkable property:
any right $D(L)$-module is projective as a right $L$-module.
\end{remarks}

\subsubsection*{Relative Hopf modules and total integrals}
Let $L$ be a Hopf algebra over a commutative ring
and $A$ a $L$-comodule algebra. Associated to this is
an entwining structure $(A, L,\psi)$, with
$$\psi(h\ot a)=a_{\psi}\ot h^{\psi}=a_{[0]}\ot ha_{[1]}$$
The resulting category of entwined modules is denoted
$$\Mm(\psi)_A^L=\Mm_A^L$$
and is usually called the category of relative Hopf modules.
Now recall \cite{Doi85} that an $L$-colinear map
$\varphi:\ L\to A$ is called
an integral. $\varphi$ is called a total integral if
$\varphi(1_L)=1_A$. We keep the notation introduced
in \eqref{4.5.1}, i.e.
\begin{equation}\eqlabel{4.10.1}
\begin{diagram}
\Mm_A^L&\rTo^F&\Mm_A\\
\dTo_H&&\dTo^{H_1}\\
\Mm^L&\rTo^{F_1}&\Mm_k
\end{diagram}
\end{equation}

\begin{theorem}\thlabel{4.10}
With notation as above, the following assertions are equivalent:\\
1) $F$ is $H$-separable;\\
2) $H_1\circ F$ is $H$-separable;\\
3) $H_1\circ F$ is $H$-Maschke;\\
4) $H_1\circ F$ is dual $H$-Maschke;\\
5) $H_1\circ F$ is $H$-semisimple;\\
6) there exists a map $\theta:\ L\ot L\to A$ such that
$\theta\circ \Delta_L=\eta_A\circ\varepsilon_L$ and
\begin{equation}\eqlabel{4.10.2}
\theta(h\ot k_{(1)})\ot k_{(2)}=\theta(h_{(2)}\ot k)_{[0]}\ot
h\theta(h_{(2)}\ot k)_{[1]}
\end{equation}
for all $h,k\in L$;\\
7) there exists a total integral $\varphi:\ L\to A$.
\end{theorem}

\begin{proof}
$2)~\Rightarrow~ 1)$: from \prref{2.3}.\\
$2)~\Rightarrow~ 3)$: from \prref{3.3}.\\
$3)~\Leftrightarrow~ 4)~\Leftrightarrow~ 5)$: from \prref{3.7}.\\
$1)~\Leftrightarrow~ 6)$: from \prref{4.8}.\\
$6)~\Rightarrow~ 7)$: Define $\varphi:\ L\to A$ by
$\varphi(h)=\theta(1\ot h)$ for all $h\in L$. A
straightforward computation shows
that $\varphi$ is a total integral.\\
$7)~\Rightarrow~ 6)$: Define $\theta:\ L\ot L\to A$ by
$\theta(h\ot k)=\varphi(S(h)k)$. It is easy to compute that
$\theta$ satisfies \eqref{4.10.2} and that
$\theta\circ \Delta_L=\eta_A\circ\varepsilon_L$.\\
$3)~\Rightarrow~ 7)$: if $H_1F$ is $H$-Maschke, then
$$H(\eta_A):\ H(A)=A\to HGK_1H_1F(A)=\Hom(A,A)\ot L$$
has a left inverse $\nu$ in $\Mm^L$, by \thref{3.4}. Now let
$\varphi(h)=\nu(I_A\ot h)$, for all $h\in L$.
Then $\varphi$ is an integral, since $\nu\in \Mm^L$, and
$$\varphi(1)=\nu(I_A\ot 1)=\nu(\eta_A(1_A))=1_A$$
$7)~\Rightarrow~ 2)$: let $\varphi$ be a total integral. We define
a natural transformation $\nu:\ HGK_1H_1F\to H$ as follows:
$$\nu_M:\ \Hom(A,M)\ot L\to M~~;~~
\nu_M(f\ot h)=f(1_A)_{[0]}\varphi(S(f(1_A)_{[1]})h)$$
We leave it to the reader to verify that $f$ is natural. Finally
\begin{eqnarray*}
\nu_M(H(\eta_M)(m))&=&
\nu_M(m_{[0]}\bullet\ot m_{[1]})\\
&=& m_{[0]}\varphi(S(m_{[1]})m_{[2]})\\
&=& m\varphi(1_H)=m
\end{eqnarray*}
\end{proof}

\subsubsection*{Doi's $[L,C]$-modules and augmented cointegrals}
Let us now discuss the dual situation. Let $L$ be a Hopf algebra,
and $C$ a right $L$-module coalgebra. We then have an entwining
structure $(L,C,\alpha)$, with
$$\alpha(c\ot h)=h_{\alpha}\ot c^{\alpha}=h_{(1)}\ot ch_{(2)}$$
The associated entwining modules are called $[L,C]$-modules, and
our diagram of forgetful functors now takes the form:
\begin{equation}\eqlabel{4.11.1}
\begin{diagram}
\Mm_L^C&\rTo^F&\Mm_L\\
\dTo_H&&\dTo^{H_1}\\
\Mm^C&\rTo^{F_1}&\Mm_k
\end{diagram}
\end{equation}
$H_1\circ F$ has a right adjoint, and this has been used in the proof
of \thref{4.10}. $H$ has a left adjoint, but, in general, $F_1$ has
no left adjoint, and this is the reason why the proof of \thref{4.11}
is different from the one of \thref{4.10}.
We recall \cite{Doi1983} that
a right $L$-linear map $\psi :C\to L$ is called
a cointegral. Furhermore, $\psi$ is called
augmented if
$\varepsilon_L \circ \psi =\varepsilon_C$.

\begin{theorem}\thlabel{4.11}
With notation as above, the following assertions are equivalent:\\
1) $H$ is $F$-separable;\\
2) $F_1\circ H$ is $F$-separable;\\
3) $F_1\circ H$ is $F$-Maschke;\\
4) $F_1\circ H$ is dual $F$-Maschke;\\
5) $F_1\circ H$ is $F$-semisimple;\\
6) there exists a map $e:\ C\to L\ot L$ such that
$m_L\circ e=\eta_L \circ \varepsilon_C$ and
\begin{equation}\eqlabel{4.11.2}
\sum e^1(c)\ot e^2(c)h=h_{(1)}e^1(ch_{(2)})\ot e^2(ch_{(2)})
\end{equation}
for all $h\in L$ and $c\in C$;\\
7) there exists an augmented cointegral $\psi:\ C\to L$.
\end{theorem}

\begin{proof}
$2)~\Rightarrow~ 1)$: from \prref{2.3};\\
$2)~\Rightarrow~ 3)$: from \prref{3.3}.\\
$3)~\Leftrightarrow~ 4)~\Leftrightarrow~ 5)$: from \prref{3.7}.\\
$1)~\Leftrightarrow~ 6)$: from \prref{4.9}.\\
$6)~\Rightarrow~ 7)$: Let $e:\ C\to L\ot L$ be as in
6); then $\psi:\ C\to L$, $\psi(c)=\varepsilon_L(e^1(c))e^2(c)$
is an augmented cointegral. \\
$7)~\Rightarrow~ 2)$: Let $\psi: C\to L$ be an augmented cointegral.
We define a natural transformation
$\Pp:\ \Hom(F_1H,F_1H)\to \Hom_L(F,F)$ as follows: for a $k$-linear
map $f:\ M\to N$, with $M,N\in \Mm_L^C$, we define
$\Pp_{M,N}(f)$ by
$$\Pp_{M,N}(f)(m)=f\Bigl(m_{[0]}S\bigl(\psi(m_{[1]})_{(1)}\bigr)\Bigr)
\psi(m_{[1]})_{(2)}$$
It is straightforward to verify that
$\Ff_1\circ\Hh=\Pp\circ \Ff$.\\
$4)~\Rightarrow~ 7)$: Assume that $F_1\circ H$ is dual $F$-Maschke.
The map
$$f:\ C\ot L\to C,~~f(c\ot h)=ch$$
is $L$-linear and $C$-colinear. As a $k$-linear map, $f$ is cosplit
by the map $c\to c\ot 1_L$, so there exists an $L$-linear map
$g:\ C\to C\ot L$, such that the following diagram commutes:
$$\begin{diagram}
C\ot L&\rTo^{f}&C\\
&\NW_g&\uTo^{I_C}\\
&&C\end{diagram}$$
Now consider $\psi=(\varepsilon_C\ot I_L)\circ g:\ C\to L$.
$\psi$ is $L$-linear. For a fixed $c\in C$, write $g(c)=\sum_i c_i\ot h_i$.
Then $c=g(f(c))=\sum_i c_ih_i$, and
\begin{eqnarray*}
\varepsilon_L(\psi(c))&=& (\varepsilon_C\ot \varepsilon_H)(g(c))\\
&=& \sum_i \varepsilon_C(c_i)\varepsilon_H(h_i)=\varepsilon_C(c)
\end{eqnarray*}
and $\psi$ is an augmented cointegral.
\end{proof}


\begin{thebibliography}{99}

\bibitem{Brzezinski1}
T. Brzezi\'nski, On modules associated to coalgebra-Galois extensions,
{\sl J. Algebra}, {\bf 215} (1999), 290--317.

\bibitem{BrzezinskiCMZ00}
T. Brzezi\'nski, S. Caenepeel, G. Militaru, and Shenglin Zhu,
Frobenius and Maschke type Theorems for Doi-Hopf modules
and entwined modules revisited: a unified approach, in
``Ring theory and Algebraic Geometry",
A. Granja, J. Hermida Alonso, and A. Verschoren (eds.),
Lecture Notes Pure Appl. Math. 221, Marcel Dekker, New York,
2001, 1-31.

\bibitem{BrzezinskiM}
T. Brzezi\'nski and S. Majid, Coalgebra bundles,
{\sl Comm. Math. Phys.} {\bf 191} (1998), 467--492.

\bibitem{book}
S. Caenepeel, G. Militaru and S. Zhu,
Frobenius and separable functors for
generalized module categories and nonlinear equations, monograph, to appear.

\bibitem{CaenepeelIMZ98}
S. Caenepeel, B. Ion, G. Militaru, and Shenglin Zhu,
Separable functors for the category of Doi-Hopf
modules, Applications, {\sl Adv. Math.} {\bf 145} (1999), 239--290.

\bibitem{CaenepeelK01}
S. Caenepeel, L. Kadison, Are Biseparable Extensions Frobenius?,
{\sl K-theory}, in press.

\bibitem{CaenepeelMZ97b}
S. Caenepeel, G. Militaru, and S. Zhu, Doi-Hopf modules,
Yetter-Drinfeld modules
and Frobenius type properties, {\sl Trans. Amer. Math. Soc.}
{\bf 349} (1997), 4311--4342.

\bibitem{CaenepeelMZ97a}
S. Caenepeel, G. Militaru and S. Zhu, A Maschke type theorem for Doi-Hopf
modules, {\sl J. Algebra} {\bf 187} (1997), 388--412.

\bibitem{Cohn}
P.M. Cohn, ``Algebra", Wiley, 1989.

\bibitem{DNT}
S. Dascalescu, C. Nastasescu and B. Torrecillas, Co-Frobenius Hopf
algebras: integrals, Doi-Koppinen modules and injective objects,
{\sl J. Algebra} {\bf 220} (1999), 542-560.

\bibitem{DI71}
F. DeMeyer and E. Ingraham, ``Separable algebras over commutative rings'',
{\sl Lecture Notes in Math.} {\bf 181}, Springer Verlag, Berlin, 1971.

\bibitem{Doi1983}
Y. Doi, On the structure of relative Hopf modules,
{\sl Comm. in Algebra} {\bf 11} (1983), 243-255.

\bibitem{Doi85}
Y. Doi, Algebras with total integral, {\sl Comm. in Algebras}, {\bf 13}
(1985),
2137--2159.

\bibitem{Kadison99b}
L. Kadison, New examples of Frobenius extensions, University Lecture
Series {\bf 14} (1999), Amer. Math. Soc., Providence.

\bibitem{KO}
M. A. Knus and M. Ojanguren, ''Th\'eorie de la descente et alg\`ebres
d'Azumaya'',
{\sl Lecture Notes in Math.} {\bf 389}, Springer Verlag, Berlin, 1974.

\bibitem{MeniniMilitaru00a}
C. Menini and G. Militaru, Integrals, quantum Galois extensions
and the affineness criterion for quantum Yetter-Drinfeld modules,
{\sl J. Algebra}, to appear.

\bibitem{Nastasescu}
C. N\u ast\u asescu, Some constructions over graded rings.
Applications, {\sl J. Algebra} {\bf 120} (1989), 119--138.

\bibitem{NastasescuRV90}
C. N\v ast\v asescu, \c S. Raianu and F. Van Oystaeyen, Modules
graded by $G$-sets, {\sl Math. Z.} {\bf 203} (1990), 605--627.

\bibitem{NastasescuVV89}
C. N\u ast\u asescu, M. Van den Bergh and F. Van Oystaeyen,
Separable functors applied to graded rings, {\sl J. Algebra}
{\bf 123} (1989), 397--413.

\bibitem{Rafael90}
M. D. Rafael, Separable functors revisited,
{\sl Comm. in Algebra} {\bf 18} (1990), 1445--1459.

\bibitem{Sweedler69}
M. E. Sweedler, ``Hopf algebras'', Benjamin, New York, 1969.

\bibitem{Takeuchi77}
M. Takeuchi, Morita Theorems for categories of comodules,
{\sl J. Fac. Sci. Univ. Tokyo} {\bf 24} (1977), 629--644.

\bibitem{Y}
D.N. Yetter, Quantum groups and representations of monoidal categories,
{\sl Math. Proc. Cambridge Philos. Soc.} {\bf 108} (1990), 261-290.

\end{thebibliography}
\end{document}